\documentclass[preprint,12pt]{elsarticle}
\usepackage{amssymb}
\usepackage{amsmath}

\usepackage{physics}
\usepackage{mathtools}

\journal{Journal name}

\begin{document}

\begin{frontmatter}

\title{Linearly scalable fast direct solver based on proxy surface method for two-dimensional elastic wave scattering by cavity}

\author[label1]{Yasuhiro Matsumoto\corref{cor1}}
\ead{matsumoto@cii.isct.ac.jp}
\cortext[cor1]{Corresponding author}
\affiliation[label1]{organization={Center for Information Infrastructure, Institute of Science Tokyo},
            addressline={2-12-1-I8-21, Ookayama, Meguro-ku},
            city={Tokyo},
            postcode={152-8550},
            country={Japan}}

\author[label2]{Taizo Maruyama}
\affiliation[label2]{organization={
    Department of Civil and Environmental Engineering, Institute of Science Tokyo},
            addressline={2-12-1 W8-W302, Ookayama, Meguro-ku},
            city={Tokyo},
            postcode={152-8550},
            country={Japan}}

\begin{abstract}
This paper proposes an $O(N)$ fast direct solver for two-dimensional elastic wave scattering problems.
The proxy surface method is extended to elastodynamics to obtain shared coefficients for low-rank approximations from discretized integral operators.
The proposed method is a variant of the Martinsson--Rokhlin--type fast direct solver.
Our variant avoids the explicit computation of the inverse of the coefficient matrix, thereby reducing the required number of matrix--matrix multiplications.
Numerical experiments demonstrate that the proposed solver has a complexity of $O(N)$ in the low-frequency range and has a highly parallel computation efficiency with a strong scaling efficiency of 70\%.
Furthermore, multiple right-hand sides can be solved efficiently; specifically, when solving problems with 180 right-hand side vectors, the processing time per vector from the second vector onward was approximately 28,900 times faster than that for the first vector.
This is a key advantage of fast direct methods.
\end{abstract}


\end{frontmatter}

\section{Introduction}
Fast direct boundary element methods (BEMs), a type of accelerated BEM,
have attracted attention \cite{preuss2022recent} since the work of Pals \cite{pals2004multipole} and Martinsson and Rokhlin \cite{martinsson2005fast}.
The BEM is a promising numerical method for wave scattering problems that include an unbounded domain.
For example, it has been used for acoustic and electromagnetic problems \cite{kirkup2019boundary, chew2001fast}.
Unlike other numerical methods such as the finite difference method or the finite element method, the BEM does not require the approximation of an infinite computational region to a finite one.
However, the conventional BEM is unsuitable for large-scale problems due to its high time complexity, which is caused by the dense coefficient matrix of the resulting system of discretized linear equations.
If the resulting system is solved using conventional direct methods such as Gaussian elimination or lower-upper (LU) decomposition, the time complexity will be $O(N^3)$, where $N$ is the number of unknowns in the system.
Fast direct BEMs are able to solve the system with $O(N)$ time complexity in the best case by means of the low-rank approximability of the off-diagonal parts of the coefficient matrix.

The combination of the fast multipole method and the iterative Krylov solver is a well--known fast BEM \cite{rokhlin1985rapid}.
Fast direct BEMs, however, are fast even in cases where the fast multipole accelerated BEM would be slow, that is, for problems whose resulting system has a relatively ill-conditioned coefficient matrix and multiple right-hand sides \cite{martinsson2005fast}.
Fast direct BEMs have been developed for two-dimensional (2D) Laplace boundary value problems (BVPs) \cite{martinsson2005fast}, 2D Helmholtz BVPs \cite{MARTINSSON2007288, matsumotoJascome2016}, 2D Helmholtz periodic BVPs \cite{gillman2013fast, greengard2014fast, matsumotoJascome2019}, three-dimensional Helmholtz BVPs \cite{Greengard_Gueyffier_Martinsson_Rokhlin_2009}, and three-dimensional electromagnetic scattering \cite{electronics11223753}.

However, to the best of our knowledge, few studies (e.g., \cite{chaillat2017theory}, \cite{bagur2022improvement}) have proposed fast direct BEMs for elastic wave scattering problems.
Naive $\mathcal{H}$-matrix LU decomposition \cite{borm2003introduction} based on matrix valued adaptive cross approximation \cite{rjasanow2017matrix} has been proposed.
The time complexity of $\mathcal{H}$-matrix LU decomposition is $O({N(\log N)}^{2})$ in the best case.
Elastic wave scattering problems are BVPs of partial differential equations related to waves propagating through solids. They are important in applications such as the analysis of seismic waves \cite{BOUCHON2007157} and non-destructive testing using ultrasonic waves \cite{MARUYAMA2024426}.
In these analyses, the structures are large relative to the wavelength and linear equations with large degrees of freedom are dealt with; hence, a fast solver with an $O(N)$ time complexity is desirable.

There are two main groups of fast direct solvers capable of achieving an $O(N)$ time complexity under optimal conditions.
The first group is proxy-surface-method-based direct solvers \cite{martinsson2005fast}, which are one of the most promising types of fast direct BEM.
The second group is ULV-decomposition-based \cite{stewart1992updating} solvers.
Both approaches utilize a tree-structured decomposition of the computational domain, linear algebra, and potential theory, and are closely related to the $\mathcal{H}$-matrix framework.
However, the application of these $O(N)$ fast direct solvers to elastodynamics is non-trivial because the wave field is a vector field.
The key to achieving $O(N)$ is sharing coefficients (or bases) for the low-rank approximation of the off-diagonal parts of the $\mathcal{H}$-matrix.
To extract the shared coefficients from the discretized integral operator of elastodynamics, we adopt the proxy surface method, which has been successfully applied to the Helmholtz transmission problem \cite{matsumotoJascome2019}.
In \cite{matsumotoJascome2019}, the shared coefficients of low-rank approximations are computed separately for the solution vectors and its normal derivatives.
In elastodynamics, we compute the shared coefficients for each direction of the vector wave field.
Similar extension methods may also be applicable to ULV-based solvers.

Previous studies (\cite{chaillat2017theory} and \cite{bagur2022improvement}) have used a boundary integral equation that has fictitious eigenfrequencies at real frequencies.
Fictitious eigenfrequencies have no physical meaning and arise from the boundary integral equation formulation.
Near these eigenfrequencies, numerical results are unreliable \cite{chen1998fictitious}.
Moreover, it is impossible to know in advance exactly where the fictitious eigenfrequencies are in general cases.
To address this issue, the Burton--Miller integral equation can be used and the integrals for hypersingular kernels that need to be calculated can be dealt with using regularization and the Galerkin method \cite{yoshida2001application}.

To achieve both speed and accuracy, 
this paper proposes a fast direct BEM based on the proxy surface method formulated with a Galerkin-discretized Burton--Miller--type boundary integral equation.
As a solver to be combined with the proxy surface method, we formulate a variant of the original Martinsson--Rokhlin solver \cite{martinsson2005fast}.
The proposed variant is similar to the formula that Gillman {\it et al.} \cite{gillman2012direct} used in the proof process, and is expected to be slightly faster to compute than the original one, as it requires fewer matrix--matrix multiplications.
Moreover, since the Martinsson--Rokhlin--type fast direct solver is often used with Nystr\"om discretization, there has been little research on how to realize the algorithm for the Galerkin method.
This paper also describes in detail how the variant of Martinsson--Rokhlin--type fast direct solvers can be applied to the Galerkin method.
For simplicity, we focus on in-plane elastic wave scattering by a cavity.

The rest of this paper is organized as follows.
Section \ref{sec:statement} states the 2D elastic wave scattering problem.
Section \ref{sec:bie} shows the corresponding boundary integral equation.
Section \ref{sec:galerkin} describes the discretization of the boundary integral equation using the Galerkin method
and Section \ref{sec:preparation} introduces binary tree decomposition and blocked discretized linear equations.
Section \ref{sec:fds} and Section \ref{sec:proxy} explain the details of the proposed fast direct solver.
Section \ref{sec:numerical} presents 
numerical examples that show that the proposed method is free of fictitious eigenfrequencies (on the real axis) not only for the integral equation but also for the fast direct solver and that it has a complexity of $O(N)$.
Section \ref{sec:conc} gives the conclusions.

\section{Two-dimensional elastic wave scattering problem} \label{sec:statement}
\begin{figure}[tb]
  \centering
  \includegraphics[width=0.4\linewidth]{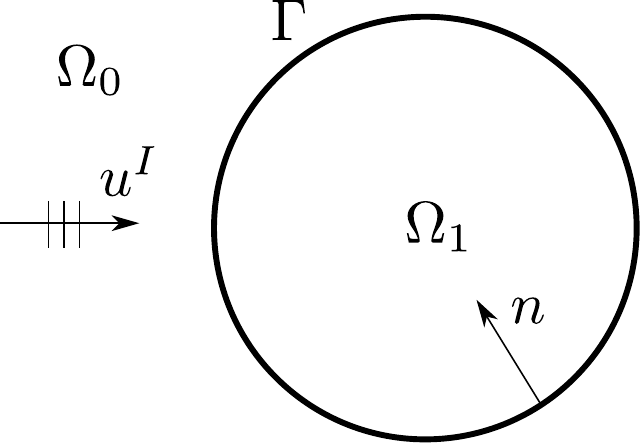}
  \caption{Diagram of elastic wave scattering by cavity.}
  \label{fig:domain}
\end{figure}

Unless otherwise specified, the Cartesian coordinate system $(x_1, x_2)$ is used to express the wave field and the summation convention is assumed for the index subscripts.
We consider the time-harmonic in-plane elastic wave scattering problem with time factor $\exp(\mathrm{i} \omega t)$, where $\mathrm{i}$, $\omega$, and $t$ are the imaginary unit, angular frequency, and time, respectively.
Consider a bounded cavity $\Omega_1 \subset \mathbb{R}^2$ whose boundary $\Gamma = \partial \Omega_1$ is a smooth, non-self-intersecting closed curve.
Let $\Omega_0 = \mathbb{R}^2 \setminus \overline{\Omega_{1}}$ be a domain filled by a linearly elastic, homogeneous, and isotropic solid with density $\rho$ and Lam\'e constants $\lambda$ and $\mu$.
In $\Omega_0$, there is an incident wave $u^{I}(x) = (u_{1}^{I}(x), u_{2}^{I}(x)) \in \mathbb{C}^2$, which is the vector field. 
Let $u(x) = (u_1 (x), u_2 (x)) \in \mathbb{C}^2$ be the total displacement vector field and $u^{S}_{i}(x) = u_{i}(x) - u^{I}_{i}(x)$ be the scattered wave field.
Let $n(x) = (n_{1}(x), n_{2}(x)) \in \mathbb{R}^2$ be a unit normal vector on $\Gamma$ toward $\Omega_1$.
The above setting is shown in Figure \ref{fig:domain}.
The problem is to find $u(x)$ that satisfies the Navier--Cauchy equation in $\Omega_{0}$ with the Neumann boundary condition for $u(x)$ and the radiation condition of elastodynamics for $u^{S}(x)$ \cite{eringen1975elastodynamics}:
\begin{align}
  &\mu u_{i,jj}(x) + (\lambda + \mu) u_{j, ji}(x) + \rho \omega^{2} u_{i}(x) = 0, &x \in \Omega_{0}, \label{bvp1} \\
  &t_{i}(x) :=  C_{ipjq} u_{j, q}(x) n_{p}(x) = T_{ij}^{n_x} u_{j}(x) = 0, &x \in \Gamma, \label{bvp2} \\
  & \lim_{|x| \to \infty} |x|^{\frac{1}{2}} \qty( \pdv{u^{S}_{i;L}(x)}{|x|} - ik_{L} u^{S}_{i;L} (x)) = 0, \\
  & \lim_{|x| \to \infty} |x|^{\frac{1}{2}} \qty( \pdv{u^{S}_{i;T}(x)}{|x|} - ik_{T} u^{S}_{i;T} (x)) = 0, \\
  & \lim_{|x| \to \infty} |x|^{-\frac{1}{2}} u^{S}_{i;L} (x) = 0, \\
  & \lim_{|x| \to \infty} |x|^{-\frac{1}{2}} u^{S}_{i;T} (x) = 0, \label{bvp3}
\end{align}
where $t_{i}(x)$ is the traction in the $x_{i}$ direction, $C_{ipjq} := \lambda \delta_{ip}\delta_{jq} + \mu (\delta_{ij}\delta_{pq} + \delta_{iq}\delta_{pj})$ is a component of the elastic tensor, $\delta_{ij}$ is Kronecker's delta, and $T_{ij}^{n_x} := C_{ipjq} n_{p} (x) \frac{\partial}{\partial x_{q}}$ is the traction operator with $n(x)$.
$u^{S}_{i;L}(x)$ and $u^{S}_{i;T}(x)$ are the longitudinal wave component and the transverse wave component of $u^{S}_{i}(x)$, respectively.
$k_L = \omega/c_L$ is the wavenumber of the longitudinal wave, where $c_L = \sqrt{(\lambda + 2\mu)/\rho}$ is the velocity of the longitudinal wave.
Similarly, $k_T = \omega/c_T$ is the wavenumber of the transverse wave, where $c_T = \sqrt{\mu/\rho}$ is the velocity of the transverse wave.

In what follows, arguments such as $u_{i}(x)$ are simplified to $u_{i}$ when they are clear from the context.

\section{Boundary integral equation} \label{sec:bie}
We now describe the boundary integral equation corresponding to the BVP \eqref{bvp1}--\eqref{bvp3}.
To avoid the fictitious eigenfrequency problem, we employ the Burton--Miller--type boundary integral equation \cite{burton1971application}:
\begin{align}
  \qty(\mathcal{D}_{ij} + \frac{\delta_{ij}}{2}  + \alpha \mathcal{N}_{ij}) u_{j} = u^{I}_{i} + \alpha T_{ij}^{n_x} u^{I}_{j} \quad \text{on } \Gamma \label{eq:bie},
\end{align}
where $\alpha \in \mathbb{C}$ is a constant, $\mathcal{D}_{ij}$ and $\mathcal{N}_{ij}$ are respectively integral operators related to the double-layer potential and a traction of the double-layer potential defined by
\begin{align}
  &\mathcal{D}_{ij} v_{j} (x) := \mathrm{v.p.} \int_{\Gamma} \qty(T_{jk}^{n_y} G_{ik} (x - y)) v_{j}(y) ds(y), \label{eq:dlp} \\
  &\mathcal{N}_{ij} v_{j} (x) := \mathrm{p.f.} \int_{\Gamma} \qty(T_{ik}^{n_x} \qty(T_{jp}^{n_y} G_{kp} (x - y))) v_{j}(y) ds(y), \label{eq:ddlp}
\end{align}
in which $G_{ij}$ is the fundamental solution of elastodynamics in two dimensions defined by
\begin{multline}
  G_{ij} (x - y) = \frac{\mathrm{i}}{4 \mu} \left( H_{0}^{(1)}(k_{T} |x - y|) \delta_{ij} \right. \\
  \left. + \frac{1}{k_{T}^{2}} \left( H_{0}^{(1)}(k_{T} |x - y|) - H_{0}^{(1)}(k_{L} |x - y|) \right)_{,ij} \right), \quad x, y \in \mathbb{R}^2. \label{eq:funda}
\end{multline}
The integrals with v.p.\@ and p.f.\@ are the principal value integral and the finite part of the divergence integral, respectively.
In \eqref{eq:funda}, $H_{0}^{(1)}$ is the Hankel function of the first kind and the zeroth order.
In \eqref{eq:bie}, the single-layer potential and its traction are eliminated by boundary condition \eqref{bvp2}.

After the solution $u_{i}$ on the boundary $\Gamma$ is obtained,
the solution of BVP \eqref{bvp1}--\eqref{bvp3} in domain $\Omega_{0}$ can be obtained using integral representation:
\begin{align}
  u_i(x) = u^{I}_{i}(x) - \int_{\Gamma} \qty(T_{jk}^{n_y} G_{ik} (x - y)) u_{j}(y) ds(y), \quad x \in \Omega_{0}. \label{eq:integ_rep}
\end{align}
In this study, we discuss the fast method up to solving the boundary integral equation \eqref{eq:bie}.
The proposed fast method is applicable to the efficient computation of \eqref{eq:integ_rep} using the method in \cite{matsumoto2024fast}.

\section{Galerkin discretization} \label{sec:galerkin}
We discretize \eqref{eq:bie} using the Galerkin method because \eqref{eq:bie} has a hypersingular kernel in the integral operator $\mathcal{N}_{ij}$.
It has been pointed out that when the collocation method is employed for crack analysis that has a hypersingular kernel in the integral operator, the accuracy of the solution will deteriorate unless boundary elements with $C^1$ class or higher smoothness are used; in contrast, when the Galerkin discretization method is employed, $C^0$ class boundary elements and test functions are acceptable \cite{yoshida2001application}.
For $C^0$ class boundary elements, the piecewise linear basis is available.
Therefore, it is expected that by combining the Burton--Miller method and the Galerkin method, the proposed solver will be easy to implement and highly accurate.

First, boundary $\Gamma$ is approximated as an $N$-polygon, where $N$ is the number of nodes. Let $\{ x^{t} \}_{t = 1}^{N}$ be a set of nodes whose coordinates are the vertices of an $N$-polygon.
Let $\{ \varphi^{t} (x) \}_{t = 1}^{N}$ be a set of piecewise linear bases on the approximated $\Gamma$.
By using $\varphi^{t}$, the unknown functions $u_i$ of \eqref{eq:bie} on the boundary can be approximated as
\begin{align}
  u_{i} (x) = \sum_{t = 1}^{N} \varphi^{t} (x) u_{i}^{t},
\end{align}
where $u^{t} = (u_{1}^{t}, u_{2}^{t}) = (u_{1} (x^t), u_{2} (x^t)) \in \mathbb{C}^2$.
Let $\{ \psi^{s} \}_{s = 1}^{N}$ be a set of test functions such that $\psi^{s} = \varphi^{s}$ for $s = 1, 2, \ldots , N$, for simplicity.
By using the $s$-th test function $\psi^{s}$ and basis functions, \eqref{eq:bie} can be formally discretized as
\begin{multline}
  \int_{\Gamma} \psi^{s}(x) \qty( \qty(\mathcal{D}_{ij} + \frac{\delta_{ij}}{2} + \alpha \mathcal{N}_{ij}) \sum_{t = 1}^{N} \varphi^{t}(y) u_{j}^{t} ) ds(x) \\
  = \int_{\Gamma} \psi_{}^{s}(x) \qty(u^{I}_{i}(x) + \alpha T_{ij}^{n_x} u^{I}_{j}(x)) ds(x), \quad x \in \Gamma, \label{eq:gal}
\end{multline}
for $s = 1, 2, \ldots, N$.
The simplest way to find the discretized solution $u$ in this Galerkin method is to solve the system given by \eqref{eq:gal} created by $N$ test functions $\{ \psi^{s} \}_{s = 1}^{N}$, which corresponds to the ordinary BEM.

\section{Preliminaries} \label{sec:preparation}
\subsection{Binary tree decomposition}
\begin{figure}[tb]
  \centering
  \includegraphics[width=0.6\linewidth]{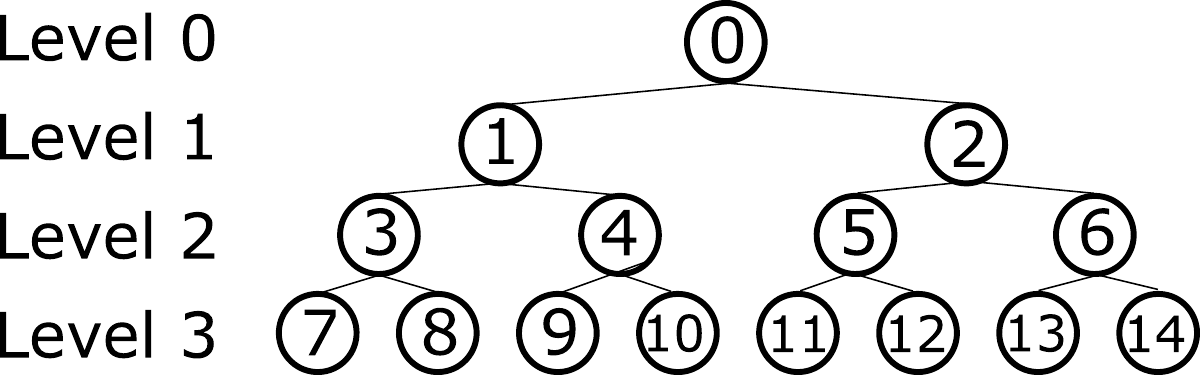}
  \caption{Binary tree decomposition with level 3 as deepest level.
    The root of the binary tree is level $0$ and has cell index $0$.
    A parent cell whose index is $i$ has two child cells whose indices are $2i+1$ and $2i+2$.
  }
  \label{fig:binary}
\end{figure}
Before describing the proposed fast direct solver using the proxy surface method, 
we introduce a binary tree for clustering the approximated $\Gamma$.
Figure \ref{fig:binary} shows a case of a level $3$ binary tree.
Let level $0$ be the root of the binary tree.
If we call the nodes of the binary tree cells, cell $0$ is divided into two cells that have indices $1$ and $2$ at level $1$.
Then, each cell at level $1$ is divided into two cells at level $2$; that is, the cell with index $1$ is divided into cells $3$ and $4$ and the cell with index $2$ is divided into cells $5$ and $6$.
Let the deepest level of the binary tree be $L$ for this repeated uniform cluster division.
The number of cells at level $L$ is $2^{L}$.
Since the sum of the number of cells from level $0$ to level $L - 1$ is $\sum_{\ell = 0}^{L - 1} 2^{\ell} = 2^{L} - 1$,
the cell index at level $L$ ranges from $2^L -1$ to $2^L - 1 + 2^L - 1 = 2^{L+1} - 2$.
When a binary tree structure is created in this way, the index of a child cell can be determined from the index of its parent cell, and vice versa.
For example, if the index of a parent cell is $i$, the indices of its child cells are $2i+1$ and $2i+2$.

For simplicity, suppose that all cells at level $L$ have $n$ elements; therefore, $N = n 2^{L}$.
Let $J_{p}$ be a finite index sequence such that 
\begin{align}
  J_{p}(s) := s + \qty(p - (2^{L} - 1) )n, \label{eq:local_to_global}
\end{align}
for $s = 1, 2, \ldots, |J_{p}|$ and $p = 2^{L} - 1, 2^{L}, \ldots, 2^{L+1} - 2$,
where the subscript $p$ is the index of cell and $|J_{p}|$ is the number of elements of $J_{p}$.
\eqref{eq:local_to_global} indicates that local index $s$ in the $p$-th cell is converted into global index $s + \qty(p - (2^{L} - 1) )n$ at level $L$.
Note that the proposed method can be performed even if $n$ is not fixed.

\subsection{Discretized linear equation}
In what follows, $x$ and $y$ are used not only for denoting coordinates, but also for denoting solutions to linear equations.
No summation convention is used for cell indices in the discretized linear equations.

Based on the binary tree decomposition scheme, we want to assemble \eqref{eq:gal} as a linear equation at level $L$ with the Galerkin method as
\begin{align}
  Ax=f, \label{eq:ax=f}
\end{align}
where $A \in \mathbb{C}^{2N \times 2N}$, $x \in \mathbb{C}^{2N}$, and $f \in \mathbb{C}^{2N}$.
Then, \eqref{eq:ax=f} is partitioned into block form as
\begin{align}
  \mqty(
  A_{2^{L}-1} & A_{(2^{L}-1)(2^{L})} & \cdots & A_{(2^{L}-1)(2^{L+1} - 2)} \\
  A_{(2^{L}) (2^{L}-1)} & A_{2^{L}} & \cdots & A_{(2^{L})(2^{L+1} - 2)} \\
  \vdots & \vdots & \cdots & \vdots \\
  A_{(2^{L+1}-2)(2^{L}-1)} & A_{(2^{L+1}-2)(2^{L})} & \cdots & A_{2^{L+1}-2} \\
  )
  \mqty(
  x_{2^L - 1} \\
  x_{2^L} \\
  \vdots \\
  x_{2^{L+1} - 2}
  )
  = 
  \mqty(
  f_{2^L - 1} \\
  f_{2^L} \\
  \vdots \\
  f_{2^{L+1} - 2}
  ), \label{eq:block_ax=f}
\end{align}
where the diagonal block part is denoted as $A_{p}$ instead of $A_{pp}$,
$A_{pq} \in \mathbb{C}^{2 |J_p| \times 2 |J_q|}$, $x_p \in \mathbb{C}^{2|J_p|}$, and $f_p \in \mathbb{C}^{2|J_p|}$.
$x_{p}$ in \eqref{eq:block_ax=f} is defined as
\begin{align}
  x_{p} :=
  \mqty(
  x_{1}^{J_{p}} \\
  x_{2}^{J_{p}}
  ), \quad
  \qty(x_{i}^{J_{p}})_{s} := u_{i}^{J_{p}(s)}
  , \label{eq:discre_x}
\end{align}
for the index of cell $p = 2^L - 1, 2^L, \ldots, 2^{L+1} - 2$ and local index $s = 1, 2, \ldots, |J_{p}|$ in the $p$-th cell.
In \eqref{eq:discre_x}, $(z)_{s}$ is the $s$-th component of vector $z$.
Similarly, $f_{p}$ in \eqref{eq:block_ax=f} is defined as
\begin{align}
  f_{p} :=
  \mqty(
  f_{1}^{J_{p}} \\
  f_{2}^{J_{p}}
  ), \quad
  \qty(f_{i}^{J_{p}})_{s} := \int_{\Gamma} \psi_{}^{J_{p}(s)}(x) \qty(u^{I}_{i}(x) + \alpha T_{ij}^{n_x} u^{I}_{j}(x)) ds(x),
\end{align}
for the index of cell $p = 2^L - 1, 2^L, \ldots, 2^{L+1} - 2$ and local index $s = 1, 2, \ldots, |J_{p}|$ in the $p$-th cell.
$A_{pq}$ in \eqref{eq:block_ax=f} is defined as
\begin{align}
  A_{pq} := D^{J_{p}J_{q}} + \frac{I^{J_{p}J_{q}}}{2} + \alpha N^{J_{p}J_{q}},
\end{align}
where $J_{p}$ and $J_{q}$ are the index sequences for test functions and basis functions with respect to the $p$-th and $q$-th cells, respectively.
The matrices in the above equation, namely $D^{J_{p}J_{q}}$, $I^{J_{p}J_{q}}$, and $N^{J_{p}J_{q}}$, are respectively defined as
\begin{align}
  D^{J_{p}J_{q}} :=
  \mqty(
  D_{11}^{J_{p}J_{q}} & D_{12}^{J_{p}J_{q}} \\
  D_{21}^{J_{p}J_{q}} & D_{22}^{J_{p}J_{q}}
  ),
\end{align}
\begin{align}
  \quad
  I^{J_{p}J_{q}} :=
  \mqty(
  I_{0}^{J_{p}J_{q}} &  \\
   & I_{0}^{J_{p}J_{q}}
  ),
\end{align}
\begin{align}
  \quad
  N^{J_{p}J_{q}} :=
  \mqty(
  N_{11}^{J_{p}J_{q}} & N_{12}^{J_{p}J_{q}} \\
  N_{21}^{J_{p}J_{q}} & N_{22}^{J_{p}J_{q}}
  ),
\end{align}
where the $(s, t)$ components of matrices $D_{ij}^{J_{p}J_{q}}$, $I_{0}^{J_{p}J_{q}}$, and $N_{ij}^{J_{p}J_{q}}$ are respectively defined as
\begin{multline}
  \qty(D_{ij}^{J_{p}J_{q}})_{st} := \int_{\Gamma} \psi_{}^{J_{p}(s)}(x) 
   \qty(\mathrm{v.p.} \int_{\Gamma} \qty( T_{jk}^{n_y} G_{ik} (x - y) ) \varphi_{}^{J_{q}(t)}(y) ds(y)) ds(x),
\end{multline}
\begin{align}
  \qty(I_{0}^{J_{p}J_{q}})_{st} := 
  \int_{\Gamma} \psi_{}^{J_{p}(s)}(x) \varphi^{J_{q}(t)}(x) ds(x),
\end{align}
\begin{multline}
  \qty(N_{ij}^{J_{p}J_{q}})_{st} := \int_{\Gamma} \psi_{}^{J_{p}(s)}(x)
   \qty(\mathrm{p.f.} \int_{\Gamma} \qty( T_{ik}^{n_x} \qty( T_{jr}^{n_y} G_{kr} (x - y) ) ) \varphi_{}^{J_{q}(t)}(y) ds(y)) ds(x),
\end{multline}
for $s = 1, 2, \ldots, |J_{p}|$ and $t = 1, 2, \ldots, |J_{q}|$.
Similarly, we define $A_{ij}^{J_{p}J_{q}}$ as
\begin{align}
  A_{ij}^{{J_{p}}{J_{q}}} := D_{ij}^{{J}_{p}{J}_{q}} + \alpha N_{ij}^{{J}_{p}{J}_{q}}. \label{eq:Aijpq}
\end{align}

To compute the discretized integral operator $\qty(N_{ij}^{J_{p}J_{q}})_{st}$, we employ a regularization technique from Appendix L.5 of \cite{YosK:2001}.
This regularization technique \cite{YosK:2001} can be easily modified for 2D cases.

\section{Fast direct solver} \label{sec:fds}
\subsection{Low-rank approximation of off-diagonal blocks}
This paper employs a modified version of the fast direct solver proposed by Martinsson and Rokhlin \cite{martinsson2005fast} for solving \eqref{eq:ax=f}.
The modified fast direct solver has the same complexity as that of the original version, namely $O(N)$ in the best case.
Unlike the original method, the modified version does not explicitly compute $A^{-1}$ in \eqref{eq:ax=f}; instead, it recursively transforms it into a linear equation with fewer degrees of freedom.
This increases computational speed because fewer matrix--matrix multiplications are required compared to those for the original method based on the Sherman--Morrison--Woodbury formula.
Although the formulas of the proposed solver are equivalent to those in the proof given by Gillman {\it et al.}\@ \cite{gillman2012direct}, the latter did not use these formulas in their numerical method;
instead, they applied the Woodbury identity, as was done in the original method.

At the deepest level $L$ of the binary tree,
from block form \eqref{eq:block_ax=f} of the linear equation \eqref{eq:ax=f},
we can consider the $i$-th block row of the linear equation as
\begin{align}
  A_{i}x_{i} + \sum_{j \neq i} A_{ij} x_{j} = f_{i}, \label{eq:block_row}
\end{align}
with respect to $i, j \in \{2^{L}-1, 2^{L}, \ldots,  2^{L+1} - 2 \}$.
Suppose that an off-diagonal block $A_{ij}$ $(i \neq j)$ can be low-rank approximated as
\begin{align}
  A_{ij} = U_{i} R_{ij} V_{j}, \label{eq:lowrank}
\end{align}
where $R_{ij}$ is the skeleton of $A_{ij}$ defined as
\begin{align}
  R_{ij} :=
  \mqty(
  A_{11}^{\tilde{J_{i}}\tilde{J_{j}}} & A_{12}^{\tilde{J_{i}}\tilde{J_{j}}} \\
  A_{21}^{\tilde{J_{i}}\tilde{J_{j}}} & A_{22}^{\tilde{J_{i}}\tilde{J_{j}}}
  )
  \in \mathbb{C}^{2 k \times 2 k}, \quad \text{for } i \neq j,
\end{align}
where a low-rank approximated rank $k$ ($k < n$) is used, $\tilde{J}_{i}$ is a sub-sequence of $J_{i}$ and $| \tilde{J}_{i} | = k$, and $A_{ij}^{\tilde{J_{p}}\tilde{J_{q}}}$ is that in \eqref{eq:Aijpq}.
In this paper, $k$ is assumed to be a fixed rank for simplicity; however, it can easily be changed to a dynamic rank.
In fact, our implementation uses a dynamic rank.
Skeleton $R_{ij}$ is a submatrix that consists of the matrix entries at the intersection of $2k$ rows and $2k$ columns of matrix $A_{ij}$.
$U_{i}$ and $V_{j}$ are coefficients of linear combinations for $R_{ij}$ defined by
\begin{align}
U_{i} = \mqty(U_{i}^{1} & \\ & U_{i}^{2}) \in \mathbb{C}^{2 |J_i| \times 2k}, \quad
V_{i} = \mqty(V_{i}^{1} & \\ & V_{i}^{2}) \in \mathbb{C}^{2 k \times 2 |J_j|}, \label{eq:uv}
\end{align}
where $U_{i}^{1} \in \mathbb{C}^{|J_i| \times k}$ and $U_{i}^{2} \in \mathbb{C}^{|J_i| \times k}$ are the coefficients of the linear combination of the first and second row block  of $R_{ij}$, respectively.
Similarly, $V_{i}^{1} \in \mathbb{C}^{k \times |J_i|}$ and $V_{i}^{2} \in \mathbb{C}^{k \times |J_i|}$ are the coefficients of the linear combination of the first and second column blocks of $R_{ij}$, respectively.
In addition, when we fix the cell index $i$ at level $\ell$, $U_i$ is shared across all $R_{ij}$ at level $\ell$ if $j \neq i$.
Similarly, for fixed $j$, $V_j$ at level $\ell$ is shared across $R_{ij}$ at level $\ell$ as well if $i \neq j$.
We explain the construction of these assumed low-rank approximations in Section \ref{sec:proxy}.

We now explain the algorithm for the fast direct solver, which has two steps.
The first step, called ``upward'', compresses the linear equation using low-rank approximations of the off-diagonal blocks and linear algebra operations.
The main part of the upward step can be performed without depending on the right-hand side.
The second step, called ``downward'', converts the solution obtained from the compressed linear equation into the original solution at level $L$.
As we will show in numerical experiments, the proposed method is suitable for problems that have multiple right-hand sides because the downward step is very fast.

\subsection{Upward step} \label{sec:upward}
Now, we convert \eqref{eq:block_row} into an equation with fewer degrees of freedom using \eqref{eq:lowrank} as follows:
\begin{align}
  A_i x_i + \sum_{j \neq i} A_{ij} x_{j} &= f_{i} \\
  A_i x_i + \sum_{j \neq i} U_{i} R_{ij} V_{j} x_{j} &= f_{i} \\
  x_i + A_{i}^{-1} \sum_{j \neq i} U_{i} R_{ij} V_{j} x_{j} &= A_{i}^{-1} f_{i} \label{eq:elim} \\
  x_i + A_{i}^{-1} U_{i} \sum_{j \neq i} R_{ij} y_{j} &= A_{i}^{-1} f_{i} &(\because y_i := V_i x_i) \\
  V_i x_i + V_i A_{i}^{-1} U_{i} \sum_{j \neq i} R_{ij} y_{j} &= V_i A_{i}^{-1} f_{i} \\
  y_i + V_i A_{i}^{-1} U_{i} \sum_{j \neq i} R_{ij} y_{j} &= V_i A_{i}^{-1} f_{i} \\
  \tilde{A}_i y_i + \sum_{j \neq i} R_{ij} y_{j} &= \tilde{A}_i V_i A_{i}^{-1} f_{i} &(\because \tilde{A}_i := (V_i A_{i}^{-1} U_{i})^{-1}) \\
  \tilde{A}_i y_i + \sum_{j \neq i} R_{ij} y_{j} &= \tilde{f}_i &(\because \tilde{f}_i := \tilde{A}_i V_i A_{i}^{-1} f_{i}) \label{eq:linear_last}
\end{align}
These row-block-wise equations have a block size of $2k$, which is smaller than $2n$ (or $2 |J_{i}|$).

By rearranging the above row-block-wise equations, the compression scheme can be recursively applied just like in the original method \cite{martinsson2005fast}.
Let $\ell$ be the level of the binary tree at which compression should be performed.
Vector-valued governing equations, such as the Navier--Cauchy equations, require care in how they are rearranged.
Let $\tilde{A}_{ij}^{p} \in \mathbb{C}^{k \times k}$ be a sub-block of $\tilde{A}_{p} \in \mathbb{C}^{2k \times 2k}$, expressed as
\begin{align}
  \tilde{A}_{p} =
  \mqty(
  \tilde{A}_{11}^{p} & \tilde{A}_{12}^{p} \\
  \tilde{A}_{12}^{p} & \tilde{A}_{22}^{p}
  ).
\end{align}
Here, the block-wise equation for level $\ell$ is rearranged into that for level $\ell - 1$ such that diagonal block $A_{i} \in \mathbb{C}^{4k \times 4k}$ is computed as
\begin{align}
  A_{i} =
  \mqty(
  \tilde{A}_{11}^{2i+1} & A_{11}^{\tilde{J}_{2i+1} \tilde{J}_{2i+2}}   & \tilde{A}_{12}^{2i+1} & A_{12}^{\tilde{J}_{2i+1} \tilde{J}_{2i+2}} \\
  A_{11}^{\tilde{J}_{2i+2} \tilde{J}_{2i+1}} & \tilde{A}_{11}^{2i+2} & A_{12}^{\tilde{J}_{2i+2} \tilde{J}_{2i+1}} & \tilde{A}_{12}^{2i+2} \\
  \tilde{A}_{21}^{2i+1} & A_{21}^{\tilde{J}_{2i+1} \tilde{J}_{2i+2}} &   \tilde{A}_{22}^{2i+1} & A_{22}^{\tilde{J}_{2i+1} \tilde{J}_{2i+2}}  \\
  A_{21}^{\tilde{J}_{2i+2} \tilde{J}_{2i+1}} & \tilde{A}_{21}^{2i+2} & A_{22}^{\tilde{J}_{2i+2} \tilde{J}_{2i+1}} & \tilde{A}_{22}^{2i+2}
  ), \label{eq:A_rearrangement}
\end{align}
for the $i$-th cell at level $\ell - 1$.
It is a rearrangement of 
\begin{align}
  \mqty(
  \tilde{A}_{2i+1} & R_{(2i+1)(2i+2)} \\
  R_{(2i+2)(2i+1)} & \tilde{A}_{2i+2}
  ).
\end{align}
Similarly, off-diagonal block $A_{ij} \in \mathbb{C}^{4k \times 4k}$ $(i \neq j)$ is computed as
\begin{align}
  A_{ij} = 
  \mqty(
  A_{11}^{\tilde{J}_{2i+1} \tilde{J}_{2j+1}} & A_{11}^{\tilde{J}_{2i+1} \tilde{J}_{2j+2}} & A_{12}^{\tilde{J}_{2i+1} \tilde{J}_{2j+1}} & A_{12}^{\tilde{J}_{2i+1} \tilde{J}_{2j+2}} \\
  A_{11}^{\tilde{J}_{2i+2} \tilde{J}_{2j+1}} & A_{11}^{\tilde{J}_{2i+2} \tilde{J}_{2j+2}} & A_{12}^{\tilde{J}_{2i+2} \tilde{J}_{2j+1}} & A_{12}^{\tilde{J}_{2i+2} \tilde{J}_{2j+2}} \\
  A_{21}^{\tilde{J}_{2i+1} \tilde{J}_{2j+1}} & A_{21}^{\tilde{J}_{2i+1} \tilde{J}_{2j+2}} & A_{22}^{\tilde{J}_{2i+1} \tilde{J}_{2j+1}} & A_{22}^{\tilde{J}_{2i+1} \tilde{J}_{2j+2}} \\
  A_{21}^{\tilde{J}_{2i+2} \tilde{J}_{2j+1}} & A_{21}^{\tilde{J}_{2i+2} \tilde{J}_{2j+2}} & A_{22}^{\tilde{J}_{2i+2} \tilde{J}_{2j+1}} & A_{22}^{\tilde{J}_{2i+2} \tilde{J}_{2j+2}} \\
  ),
\end{align}
for the $i$-th and $j$-th cells at level $\ell - 1$.
It is a rearrangement of 
\begin{align}
  \mqty(
  R_{(2i+1)(2j+1)} & R_{(2i+1)(2j+2)} \\
  R_{(2i+2)(2j+1)} & R_{(2i+2)(2j+2)}
  ).
\end{align}
For the $i$-th cell at level $\ell - 1$, $x_i$ and $f_i$ are respectively assigned as the $s$-th component as
\begin{align}
  (x_i)_{s} =
  \begin{dcases}
    (y_{2i+1})_{s} & \text{if } 1 \leq s \leq k, \\
    (y_{2i+2})_{s - k} & \text{if } k+1 \leq s \leq 2k, \\
    (y_{2i+1})_{s - k} & \text{if } 2k+1 \leq s \leq 3k, \\
    (y_{2i+2})_{s - 2k} & \text{if } 3k+1 \leq s \leq 4k,
  \end{dcases}
  \\
  (f_i)_{s} =
  \begin{dcases}
    (\tilde{f}_{2i+1})_{s} & \text{if } 1 \leq s \leq k, \\
    (\tilde{f}_{2i+2})_{s - k} & \text{if } k+1 \leq s \leq 2k, \\
    (\tilde{f}_{2i+1})_{s - k} & \text{if } 2k+1 \leq s \leq 3k, \\
    (\tilde{f}_{2i+2})_{s - 2k} & \text{if } 3k+1 \leq s \leq 4k,
  \end{dcases}
  \label{eq:rhs_rearrangement}
\end{align}
for $s = 1, 2, \ldots, 4k$.
They are rearrangements of
\begin{align}
  \mqty(
  y_{2i+1} \\
  y_{2i+2}
  ), \quad
  \mqty(
  \tilde{f}_{2i+1} \\
  \tilde{f}_{2i+2}
  ).
\end{align}
In summary, index sequence $J_{i}$ at parent level $\ell - 1 \neq L$ is made by using child cells $\tilde{J}_{2i+1}$ and $\tilde{J}_{2i+2}$ at level $\ell$ as
\begin{align}
  J_{i}(s) = 
  \begin{dcases}
    \tilde{J}_{2i+1}(s)   & \text{if } 1 \leq s \leq k \,\,(= |\tilde{J}_{2i+1}(s) |), \\
    \tilde{J}_{2i+2}(s-k) & \text{if } k + 1 \leq s \leq 2k \,\,(= |\tilde{J}_{2i+1}(s) | + |\tilde{J}_{2i+2}(s) |).
  \end{dcases}
\end{align}
for local index $s = 1, 2, \ldots, 2k \,\,(= |J_{i}|)$.
At level $\ell - 1$, off-diagonal block $A_{ij}$ of the equation
$A_{i}x_{i} + \sum_{j \neq i} A_{ij} x_{j} = f_{i}$
can also be low-rank approximated because $A_{ij}$ at level $\ell - 1$ is a submatrix of the off-diagonal matrices at level $\ell$ constructed from integrals of the fundamental solutions (for details, see \cite{martinsson2005fast}).

\subsection{Downward step} \label{sec:downward}
By recursively executing the scheme in Section \ref{sec:upward},
we can obtain the linear equation $Ax=f$ at level $\ell$ with sufficiently small degrees of freedom.
Therefore, we can solve this compressed linear equation with little computational effort.
Since there are $2^{\ell}$ cells at level $\ell$, the compressed linear equation is expressed as
\begin{align}
  \mqty(
  A_{2^{\ell} - 1} & A_{(2^{\ell} - 1)(2^{\ell})} & \cdots & A_{(2^{\ell}-1)(2^{\ell+1} -2)} \\
  A_{(2^{\ell}) (2^{\ell}-1)} & A_{2^{\ell}} & \cdots & A_{(2^{\ell})(2^{\ell+1} - 2)} \\
  \vdots & \vdots & \cdots & \vdots \\
  A_{(2^{\ell+1}-2)(2^{\ell}-1)} & A_{(2^{\ell+1}-2)(2^{\ell})} & \cdots & A_{2^{\ell+1}-2} \\
  )
  \mqty(
  x_{2^\ell -1} \\
  x_{2^\ell } \\
  \vdots \\
  x_{2^{\ell+1} - 2}
  )
  = 
  \mqty(
  f_{2^\ell - 1} \\
  f_{2^\ell} \\
  \vdots \\
  f_{2^{\ell+1} - 2}
  ). \label{eq:block_ax=f_at_ell}
\end{align}
By solving this equation, $x_i \in \mathbb{C}^{4k}$ can be obtained for $i = 2^{\ell} - 1, 2^{\ell}, \ldots, 2^{\ell + 1} - 2$.

The solved $x_i$ at level $\ell$ is rearranged into $y_{2i+1}$ and $y_{2i+2}$ at level $\ell + 1$ such that
\begin{align}
  (y_{2i+1})_{s}   & = (x_{i})_{s} \quad \text{if } 1 \leq s \leq k, \label{eq:xtoy1} \\
  (y_{2i+2})_{s-k} & = (x_{i})_{s} \quad \text{if } k+1 \leq s \leq 2k, \\
  (y_{2i+1})_{s-k}  &= (x_{i})_{s} \quad \text{if } 2k+1 \leq s \leq 3k, \\
  (y_{2i+2})_{s-2k} &= (x_{i})_{s} \quad \text{if } 3k+1 \leq s \leq 4k. \label{eq:xtoy2}
\end{align}
We now describe the algorithm that converts $y_i$ at level $\ell+1$ into $x_i$.
First, the first term on the left--hand side of \eqref{eq:elim} is transpositioned to right--hand side as
\begin{align}
  A_{i}^{-1} U_{i} \sum_{j \neq i} R_{ij} V_{j} x_{j} = -x_i + A_{i}^{-1} f_{i}. \label{eq:elim2}
\end{align}
If we multiply both sides of \eqref{eq:linear_last} by $A_{i}^{-1} U_{i}$ from the left, we obtain the following equation:
\begin{align}
  A_{i}^{-1} U_{i} \tilde{A}_i y_i + A_{i}^{-1} U_{i} \sum_{j \neq i} R_{ij} y_{j} = A_{i}^{-1} U_{i} \tilde{f}_{i}. \label{eq:elim3}
\end{align}
For eliminating $R_{ij}$, to substitute \eqref{eq:elim2} into \eqref{eq:elim3}, we obtain
\begin{align}
  A_{i}^{-1} U_{i} \tilde{A}_{i} y_i - x_i + A_{i}^{-1}f_i = A_{i}^{-1} U_{i} \tilde{f}_i
\end{align}
and simplify to
\begin{align}
  x_i = A_{i}^{-1}f_i - A_{i}^{-1} U_{i} \tilde{f}_i + A_{i}^{-1} U_{i} \tilde{A}_{i} y_{i}. \label{eq:down_comv}
\end{align}
From this equation, we can obtain $x_i$ at level $\ell$ from $y_i$.

The above algorithm can be applied again by rearranging $x_i$ at level $\ell$ into $y_{2i+1}$ and $y_{2i+2}$ at level $\ell +1$ using \eqref{eq:xtoy1}--\eqref{eq:xtoy2}.
$x_i$ at level $L$ obtained using this algorithm is the solution of the original linear equation \eqref{eq:block_ax=f}.

\section{Proxy surface method} \label{sec:proxy}
This section briefly describes how to apply the proxy surface method.
The proxy surface method is used to make the fast direct solver have a complexity of $O(N)$.
For details on the ideas behind the proxy surface method, see \cite{martinsson2005fast}.

\subsection{Computation of $V$}
\begin{figure}[tb]
  \centering
  \includegraphics[width=0.3\linewidth]{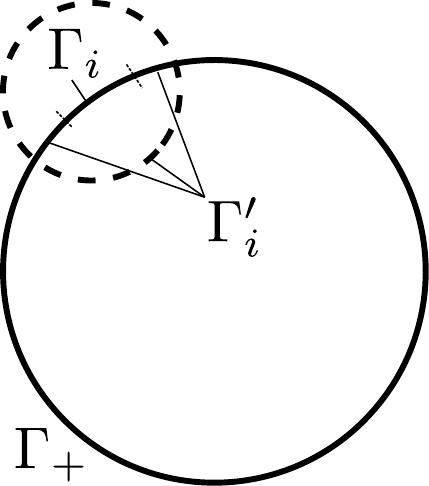}
  \caption{Proxy surface setting.
    $\Gamma_i$ is a subset of $\Gamma$ with respect to the $i$-th cell at level $\ell$.
    The proxy surface is a local virtual boundary that encloses $\Gamma_{i}$.
    $\Gamma_{i}^{\prime}$ is a union of the proxy surface and subsets of $\Gamma$ enclosed by the local virtual boundary.
    $\Gamma_{+}$ is a subset of $\Gamma$ that is outside the proxy surface.
    By using the local virtual boundary instead of $\Gamma_{+}$ to calculate the interaction with $\Gamma_i$, we can make the fast direct solver have a complexity of $O(N)$.
  }
  \label{fig:domain_proxy}
\end{figure}
Suppose that $\ell$ is the level at which the upward step is applied.
Level $\ell$ has index sequence $J_{i}$ for $i = 2^{\ell} - 1, 2^{\ell}, \ldots, 2^{\ell+1}-2$.
For a fixed $i$, let $\Gamma_{i}$ be a support of $\{ \varphi^{J_{i}(s)} \}_{s = 1}^{|J_{i}|}$ or $\{ \psi^{J_{i}(s)} \}_{s = 1}^{|J_{i}|}$; that is, $\Gamma_{i} = \bigcup_{s = 1}^{|J_i|} \operatorname{supp}(\varphi_{}^{J_{i}(s)})$ or 
$\Gamma_{i} = \bigcup_{s = 1}^{|J_i|} \operatorname{supp}(\psi_{}^{J_{i}(s)})$.
$\Gamma_{i}$ is a subset of the approximated $\Gamma$ shown in Figure \ref{fig:domain_proxy}.
Let $\Gamma_{i}^{\prime}$ be a union of a local virtual boundary that encloses $\Gamma_{i}$ and subsets of $\Gamma$ enclosed by the local virtual boundary.
The local virtual boundary is called a proxy surface.
$\Gamma_{i}^{\prime}$ is shown in Figure \ref{fig:domain_proxy}.

We now consider the discretization of $\Gamma_{i}^{\prime}$.
Suppose that a local virtual boundary is discretized into an $m^{\prime}$-polygon, where $m^{\prime} = O(1)$.
Let $m$ be the sum of $m^{\prime}$ and the number of subsets of an $N$-polygon ($\Gamma$) enclosed by the local virtual boundary.
If $m = O(1)$, we can make the fast direct solver have a complexity of $O(N)$.
On $m$ elements, we define the piecewise linear basis $(\varphi^{\prime})^{s}$ and test function $(\psi^{\prime})^{s} = (\varphi^{\prime})^{s}$ for $s = 1, 2, \ldots, m$.
Let $J_{i}^{\prime}$ be an index sequence with respect to $\Gamma_{i}^{\prime}$.
$J_{i}^{\prime}(s) = s$ for local index $s = 1, 2, \ldots, m$.
By using $J_{i}^{\prime}$ and $(\psi^{\prime})^{s}$, 
we can compute the interaction matrix $M^{V_i} \in \mathbb{C}^{m \times n}$ between $\Gamma_{i}^{\prime}$ and $\Gamma_{i}$ as follows:
\begin{align}
  M_{}^{V_i} =
  \mqty(
  A_{11}^{J_{i}^{\prime}J_{i}} & A_{12}^{J_{i}^{\prime}J_{i}} \\
  A_{21}^{J_{i}^{\prime}J_{i}} & A_{22}^{J_{i}^{\prime}J_{i}}
  ), \label{eq:mv}
\end{align}
where $A_{pq}^{J_{i}^{\prime}J_{i}}$ is computed as in \eqref{eq:Aijpq}.

The method used for computing $V_{i}$, namely $V_{i}^{1}$ and $V_{i}^{2}$ in \eqref{eq:uv}, is as follows.
We define the left-half $M_{\text{left}}^{V_i}$ and the right-half $M_{\text{right}}^{V_i}$ of \eqref{eq:mv} as
\begin{align}
  M_{\text{left}}^{V_i} =
  \mqty(
  A_{11}^{J_{i}^{\prime}J_{i}} \\
  A_{21}^{J_{i}^{\prime}J_{i}}
  ),
  \quad
  M_{\text{right}}^{V_i} =
  \mqty(
  A_{12}^{J_{i}^{\prime}J_{i}} \\
  A_{22}^{J_{i}^{\prime}J_{i}}
  ).
\end{align}
We can obtain $V_{i}^{1}$ and $V_{i}^{2}$ via the interpolative decomposition \cite{cheng2005compression} of $M_{\text{left}}^{V_i}$ and $M_{\text{right}}^{V_i}$, respectively.
Interpolative decomposition is a matrix decomposition method used in the Martinsson--Rokhlin fast direct solver \cite{martinsson2005fast}.
This scheme can be computed using column-pivoted QR decomposition.
Suppose that $P$ is a permutation matrix with size $n$, $k$ is the cut-off rank, $Q_1 \in \mathbb{C}^{m \times k}$ and $Q_2 \in \mathbb{C}^{m \times (n-k)}$ are semi-unitary matrices, $R_1 \in \mathbb{C}^{k \times k}$ and $R_3 \in \mathbb{C}^{(n-k) \times (n-k)}$ are upper triangular matrices, and $R_2 \in \mathbb{C}^{k \times (n-k)}$.
Since $M_{\text{left}}^{V_i}$ (or $M_{\text{right}}^{V_i}$) has no singularity, it can be low-rank approximated.
$M_{\text{left}}^{V_i}$ is decomposed as
\begin{align}
  M_{\text{left}}^{V_i}P &=
  \mqty(
  Q_1 & Q_2
  )
  \mqty(
  R_1 & R_2 \\
  & R_3
  ) \\
  &\approx
  Q_1
  \mqty(
  R_1 & R_2 \\
  ) &(\because \norm{R_3} \approx 0) \\
  &=
  Q_1 R_1
  \mqty(
  I & R_{1}^{-1} R_2
  ) \\
  M_{\text{left}}^{V_i} &= Q_1 R_1
  \mqty(
  I & R_{1}^{-1} R_2
  )P^{T}.
\end{align}
In the above equation, $Q_1 R_1$ is the column skeleton of $M_{\text{left}}^{V_i}P$ and
$
\mqty(
I & R_{1}^{-1} R_2
)P^{T}
$
is the linear combination coefficient $V_{i}^{1}$.
Similarly, $V_{i}^{2}$ can be obtained via the interpolative decomposition of $M_{\text{right}}^{V_i}$.

By applying this scheme for all cell indices $i$ at level $\ell$, we can obtain $V_{i}$ for $i = 2^{\ell} - 1, 2^{\ell}, \ldots 2^{\ell+1}-2$.

\subsection{Computation of $U$}
For computing $U_i$ at level $\ell$, we use almost the same method as that for $V_i$, except that we swap the roles of the observation and source bases.
The interaction matrix $M_{}^{U_i}$ between $J_{i}$ and $J_{i}^{\prime}$ is computed as
\begin{align}
  M_{}^{U_i} =
  \mqty(
  A_{11}^{J_{i}J_{i}^{\prime}} & A_{12}^{J_{i}J_{i}^{\prime}} \\
  A_{21}^{J_{i}J_{i}^{\prime}} & A_{22}^{J_{i}J_{i}^{\prime}} \\
  ).
\end{align}
We define the upper-half $M_{\text{up}}^{U_i}$ and the lower-half $M_{\text{low}}^{U_i}$ of the above matrix as
\begin{align}
  M_{\text{up}}^{U_i} =
  \mqty(
  A_{11}^{J_{i}J_{i}^{\prime}} & A_{12}^{J_{i}J_{i}^{\prime}}
  ),
  \quad
  M_{\text{low}}^{U_i} =
  \mqty(
  A_{21}^{J_{i}J_{i}^{\prime}} & A_{22}^{J_{i}J_{i}^{\prime}}
  ).
\end{align}
Since $M_{\text{up}}^{U_i}$ (or $M_{\text{low}}^{U_i}$) has no singularity,
it can be low-rank approximated.
The transpose and complex conjugate of $M_{\text{up}}^{U_i}$, namely $(M_{\text{up}}^{U_i})^{H}$, is decomposed as
\begin{align}
  (M_{\text{up}}^{U_i})^{H} P &=
  \mqty(
  Q_1 & Q_2
  )
  \mqty(
  R_1 & R_2 \\
  & R_3
  ) \\
  &\approx
  Q_1
  \mqty(
  R_1 & R_2 \\
  ) &(\because \norm{R_3} \approx 0) \\
  &=
  Q_1 R_1
  \mqty(
  I & R_{1}^{-1} R_2
  ) \\
  (M_{\text{up}}^{U_i})^{H} &= Q_1 R_1
  \mqty(
  I & R_{1}^{-1} R_2
  )P^{T} \\
  M_{\text{up}}^{U_i} &=
  \qty(
  \mqty(
  I & R_{1}^{-1} R_2
  ) P^{T}
  )^{H}
  (Q_1 R_1)^{H}.
\end{align}
In the above equation, $(Q_1 R_1)^{H}$ is the row skeleton of $P M_{\text{up}}^{U_i}$ and
$
\qty(
  \mqty(
  I & R_{1}^{-1} R_2
  ) P^{T}
  )^{H}
$
is the linear combination coefficient $U_{i}^{1}$.
Similarly, $U_{i}^{2}$ can be obtained via the interpolative decomposition of $M_{\text{low}}^{U_i}$.

By applying this scheme for all cell indices $i$ at level $\ell$, 
the low-rank approximation of off-diagonal blocks of the linear equation at level $\ell$ can be achieved with shared $U_{i}$ and $V_{i}$.

\subsection{Summary of fast direct solver algorithm} \label{sec:algorithm}
Let $\ell_{0}$ be the level at which the compressed linear equation is solved.
The algorithm of the proposed fast direct solver is summarized as follows:
\begin{enumerate}
\item Make a mesh and cluster it into a uniform binary tree whose deepest level is $L$.
\item Compute index sequence $J_{i}$ for $i = 2^{L} - 1, 2^{L}, \ldots 2^{L+1}-2$.
\item (Upward step) Loop for $\ell = L, L-1, \ldots, \ell_{0} + 1$:
  \begin{enumerate}
  \item Compute $U_{i}$, $V_{i}$, and $\tilde{J}_{i}$ for $i = 2^{\ell} - 1, 2^{\ell}, \ldots 2^{\ell+1}-2$ using the method in Section \ref{sec:proxy}.
  \item Compute $A_{i}^{-1}$, $A_{i}^{-1} U_{i}$, $\tilde{A}_i$, and $\tilde{f}_i$ for $i = 2^{\ell} - 1, 2^{\ell}, \ldots 2^{\ell+1}-2$.
  \item To rearrange the linear equation into that at level $\ell - 1$, compute $A_i$ from \eqref{eq:A_rearrangement} and $f_i$ from \eqref{eq:rhs_rearrangement} for $i = 2^{\ell-1} - 1, 2^{\ell-1}, \ldots 2^{\ell}-2$.
  \end{enumerate}
\item Solve the compressed linear equation $Ax=f$ for $x$ at level $\ell_{0}$.
\item (Downward step) Loop for $\ell = \ell_{0}, \ell_{0} + 1, \ldots, L-1$:
  \begin{enumerate}
  \item Rearrange $x_i$ at level $\ell$ into $y_{2i+1}$ and $y_{2i+2}$ at level $\ell+1$ for $i = 2^{\ell} - 1, 2^{\ell}, \ldots 2^{\ell+1} - 2$ using \eqref{eq:xtoy1}--\eqref{eq:xtoy2}.
  \item Convert $y_{2i+1}$ and $y_{2i+2}$ into $x_{2i+1}$ and $x_{2i+2}$ using \eqref{eq:down_comv}.
  \end{enumerate}
\end{enumerate}
In practice, $\ell_{0}$ can be set to $0$, $1$, or $2$.

\section{Numerical examples} \label{sec:numerical}
\begin{figure}[tb]
  \centering
  \includegraphics[width=0.5\linewidth]{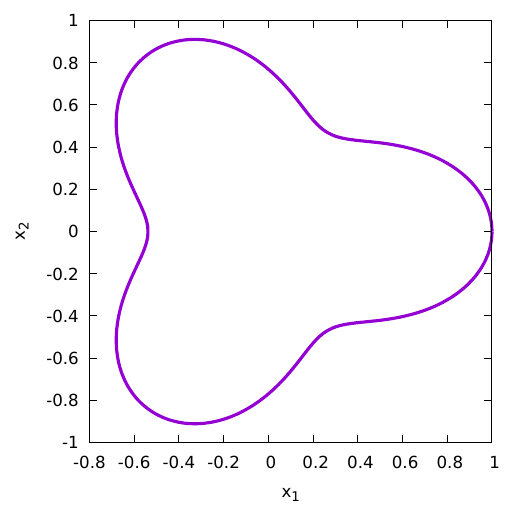}
  \caption{Scatter used in numerical examples.
    It was generated using $(x_1, x_2) = \qty((r + a \cos(b \theta)) \cos(\theta)/(1 + a), ((r + a \cos(b \theta)) \sin(\theta))/(1 + a))$ with $r = 1$, $a = 0.3$, and $b = 3$.
    It is supposed that the representative length of this mesh is 1.
  }
  \label{fig:tripod}
\end{figure}

We now test the validity and effectiveness of the proposed fast direct solver.
The computation time and accuracy of the proposed fast direct solver (FDS) are compared with those of the conventional (non-accelerated) BEM (Conv).
The algorithm shown in Section \ref{sec:algorithm} was implemented using mostly \texttt{C++} and some \texttt{Fortran}.
\texttt{Fortran} code is called from \texttt{C++} code for computing the discretized integral operators.
\texttt{Eigen} \cite{eigenweb} was used as a linear algebra library.
All numerical examples in this section use the boundary shape shown in Figure \ref{fig:tripod} as $\Gamma$.
The supercomputer Camphor3 at Kyoto University, which has two 56-core Intel Xeon Max 9480 CPUs (total number of cores: 112; maximum frequency: 1900 MHz) and 128 GB of 3.2-TB/s HBM2e memory,
was used for most calculations.
The supercomputer TSUBAME4.0 at the Institute of Science Tokyo, which has two 96-core AMD EPYC 9654 CPUs (total number of cores: 192; maximum frequency: 3708 MHz) and 768 GB of 0.92-TB/s DDR5 memory,
was used for calculations that required large amounts of memory.
A plane longitudinal wave was employed as incident wave $u^I$.
In addition, $c_{L} = \sqrt{3}$, $c_{T} = 1$, and $\rho = 1$ were used in all numerical examples.

\subsection{Relative error of proposed method} \label{sec:error}
Let $x_{\text{FDS}}$ be a solution vector computed by FDS.
Similarly, let $x_{\text{Conv}}$ be a solution vector computed by Conv.
We define the relative error as
\begin{align}
  \frac{\norm{x_{\text{FDS}} - x_{\text{Conv}}}_{2}}{\norm{x_{\text{Conv}}}_{2}},
\end{align}
where $\norm{\cdot}_{2}$ means the $p$-norm with $p = 2$.
The truncation error of the interpolative decomposition $\varepsilon$ was set to $\varepsilon = 10^{-8}$, $10^{-9}$, or $10^{-10}$.
Since our implementation of FDS uses a dynamic rank, this parameter controls the rank used for the low-rank approximation.

\begin{table}[h]
\centering
\caption{Relative error ($\norm{x_{\text{FDS}} - x_{\text{Conv}}}_{2}/{\norm{x_{\text{Conv}}}_{2}}$) of multi-level FDS for various truncation errors $\varepsilon$.
The accuracy of FDS improves as $\varepsilon$ decreases.
The decrease in accuracy with increasing $N$ is explained in the text.
}
\label{tab:error}
\begin{tabular}{llllll} \hline
$N$    & $L$& $\ell_0$& $\varepsilon=10^{-8}$              & $\varepsilon=10^{-9}$              & $\varepsilon=10^{-10}$                 \\ \hline
$400$  & 2& 1& $8.479 \times 10^{-9}$ & $4.131 \times 10^{-10}$  & $2.790 \times 10^{-11}$   \\
$800$  & 3& 1& $1.881 \times 10^{-8}$ & $2.590 \times 10^{-9}$  & $9.346 \times 10^{-11}$    \\
$1600$ & 4& 1& $5.784 \times 10^{-8}$ & $7.717 \times 10^{-9}$  & $6.191 \times 10^{-10}$    \\
$3200$ & 5& 1& $3.308 \times 10^{-7}$ & $1.801 \times 10^{-8}$  & $2.016 \times 10^{-9}$    \\
$6400$ & 6& 1& $8.404 \times 10^{-7}$ & $1.380 \times 10^{-7}$  & $1.037 \times 10^{-8}$  \\ \hline
\end{tabular}
\end{table}
Table \ref{tab:error} shows the relative errors of FDS with $N = 400, 800, 1600, 3200, 6400$.
The $N$ values were determined by fixing the cell size $n$ at 100 and setting the deepest level of the binary tree to 2, 3, 4, 5, and 6, respectively.
For example, $n2^{L} = 100 \times 2^{2} = 400$ for $n = 100$ and $L$ = 2.
Note that $N$ refers to the number of boundary divisions; therefore, the linear equation had $2N$ degrees of freedom.
The compressed linear equation $\ell_{0}$ was solved at level 1.
Dimensionless angular frequency $\omega = 2$.
Table \ref{tab:error} shows that the accuracy of FDS improves as $\varepsilon$ decreases.
These results confirm the accuracy of the proposed FDS.

However, as $N$ increases, the accuracy of FDS decreases.
There are several possible causes.
The first possible cause is the fixed cell size $n$.
There are two ways to determine the cell size.
One is to use a fixed size, as done in the current method, and
the other is to dynamically determine the cell size according to the distance between cells to satisfy low-rank approximability.
However, increasing the cell size does not prevent the accuracy of FDS from decreasing as $N$ increases.

Another possible cause is that the recursive algorithm causes the skeletons to cluster near cell boundaries.
When using the proxy surface method, skeletons near the cell boundaries are more likely to be selected because the skeletons are automatically selected via interpolative decomposition, where the column (or row) with the largest norm is selected among the columns (or rows) of the interaction matrix $M^{V_i}$ (or $M^{U_i}$) (see Section \ref{sec:proxy}).
To verify this, we executed a single-level version of FDS and measured the relative error of the solution vector; the results are shown in Table \ref{tab:error_single}.
As shown, there is almost no deterioration in the accuracy of the solution of FDS as $N$ increases.
\begin{table}[h]
\centering
\caption{Relative error ($\norm{x_{\text{FDS}} - x_{\text{Conv}}}_{2}/{\norm{x_{\text{Conv}}}_{2}}$) of single-level FDS for $\varepsilon = 10^{-8}$.
The decrease in accuracy of FDS with increasing $N$ is suppressed compared to that for the multi-level version (Table \ref{tab:error}).
}
\label{tab:error_single}
\begin{tabular}{llll} \hline
$N$    & $L$& $\ell_0$& $\varepsilon=10^{-8}$ \\ \hline
  $400$ & 2& 1& $8.479\times 10^{-9}$ \\
  $800$ & 3& 2& $2.356\times 10^{-9}$ \\
 $1600$ & 4& 3& $2.952\times 10^{-9}$ \\
 $3200$ & 5& 4& $3.348\times 10^{-9}$ \\
 $6400$ & 6& 5& $4.590\times 10^{-9}$ \\
$12800$ & 7& 6& $1.080\times 10^{-8}$ \\
$25600$ & 8& 7& $8.860\times 10^{-9}$ \\
$51200$ & 9& 8& $3.129\times 10^{-8}$ \\ \hline
\end{tabular}
\end{table}

The development of a multi-level version of FDS that is more accurate is left as a future challenge.

\subsection{$N$ versus elapsed time}
We now demonstrate the complexity of the proposed FDS.
A dimensionless angular frequency $\omega$ was set to $\omega = 2$.
The truncation error of the interpolative decomposition $\varepsilon$ was set to $\varepsilon = 10^{-8}$.
The results for $\varepsilon = 10^{-9}$ and $10^{-10}$ are omitted here since the qualitative results were unchanged (the calculation time slightly increased).

\begin{figure}[tb]
  \centering
    \includegraphics[width=0.95\linewidth]{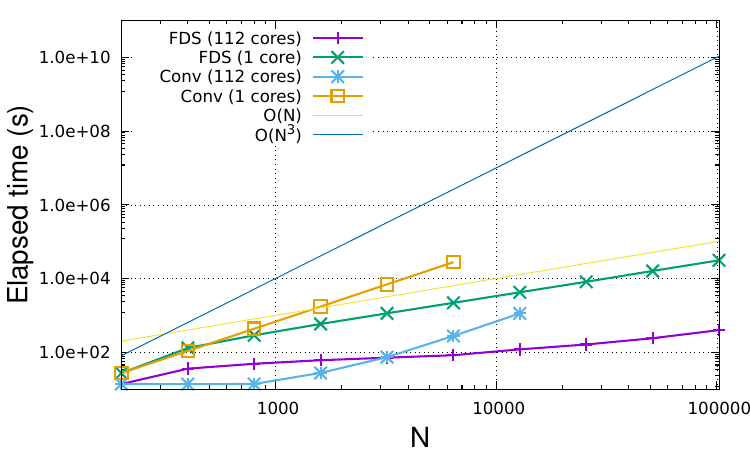}
  \caption{$N$ versus elapsed time for proposed fast direct solver (FDS) and conventional BEM (Conv).
    FDS and Conv were parallelized using \texttt{OpenMP}. The plots show results for 1 and 112 CPU cores.
    The plot of FDS (1 core) indicates that the proposed FDS has a complexity of $O(N)$.
    FDS (112 cores) is about 9.5 times faster than Conv (112 cores) at $N = 10240$.
    FDS (112 cores) is about 78 times faster than FDS (1 cores) at $N = 102400$.
    FDS (1 core) is faster than Conv (1 core) at $N \geq 800$.
    FDS (112 cores) is faster than Conv (112 cores) at $N \geq 3200$.
  }
  \label{fig:NvsTime}
\end{figure}

Figure \ref{fig:NvsTime} shows the elapsed time (in seconds) of FDS for $N = 200$ to $102400$ and that of Conv for $N = 200$ to $12800$.
Although $\ell_{0}$ and $n$ were fixed at 1 and 100, respectively, the deepest level $L$ varied with $N$.
FDS and Conv were parallelized using \texttt{OpenMP}.
FDS (1 core) shows that the proposed fast direct solver has a complexity of $O(N)$.
FDS (112 cores) seems to have a complexity of better than $O(N)$; however, this is due to the effect of parallelization since the best case complexity is $O(N)$ for the proposed FDS.
FDS also shows good strong scaling efficiency; when N=102400, it is about 78 times faster when using 112 CPU cores (i.e., strong scaling efficiency is about 70\%).
Furthermore, for $N \geq 800$, FDS (1 core) is faster than Conv (1 core).
If the fast multipole method is employed as a fast solver, a larger $N$ may be required to exceed the computational speed of Conv.

\subsection{Case with multiple right-hand sides}
We now demonstrate problems with multiple right-hand sides, for which the proposed FDS is most effective.
180 incident waves $u^I$ with different incident angles were used as multiple right-hand sides.
Dimensionless angular frequency $\omega$ was set to $\omega = 2$.
$N$ was fixed at 12800 $(L = 7)$.
In the following numerical examples, the truncation error of the interpolative decomposition $\varepsilon = 10^{-8}$ and $\ell_{0} = 1$ were fixed.

\begin{figure}[tb]
  \centering
  \includegraphics[width=0.9\linewidth]{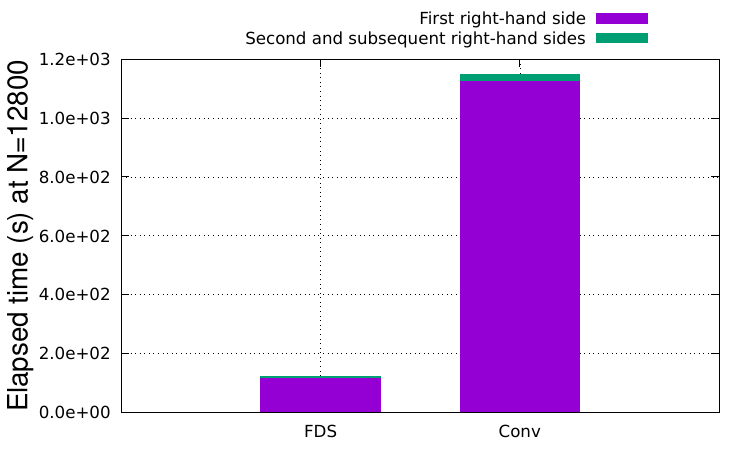}
  \caption{Elapsed time for solving 180-right-hand-side vector.
    FDS is the proposed fast direct solver.
    Conv is the conventional BEM with no fast method.
    FDS and Conv were parallelized with 112 CPU cores.
    The violet part indicates the elapsed time for solving first right-hand side.
    The green part indicates the elapsed time for solving the second and subsequent right-hand sides.
    For the second and subsequent right-hand sides, Conv seems efficient because LU decomposition is used as the solver.
  }
  \label{fig:multiple_rhs}
\end{figure}
Figure \ref{fig:multiple_rhs} shows that FDS is very efficient in solving the second and subsequent right-hand sides.
It took only 0.74 seconds to solve 179 vectors (second and subsequent right-hand sides).
This is less than 1\% of the time it took to solve the first vector (119.4 seconds).
The processing time per vector from the second vector onward was approximately 28,900 times faster than that for the first vector.

\subsection{Effect of angular frequency on performance}
In this section, we demonstrate the performance of the proposed FDS for various values of dimensionless angular frequency $\omega$.
\begin{figure}[tb]
  \centering
  \includegraphics[width=1.0\linewidth]{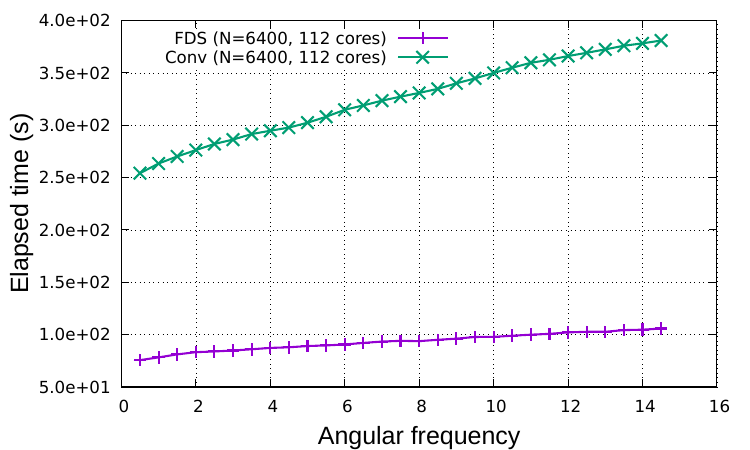}
  \caption{Elapsed time for FDS (proposed fast direct solver) and Conv (conventional BEM) for various values of angular frequency $\omega$.
  }
  \label{fig:omega_time}
\end{figure}
\begin{figure}[tb]
  \centering
  \includegraphics[width=1.0\linewidth]{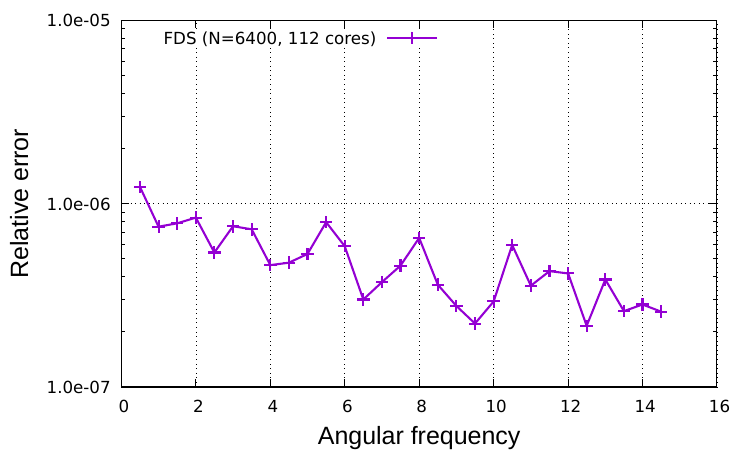}
  \caption{Relative error ($\norm{x_{\text{FDS}} - x_{\text{Conv}}}_{2}/{\norm{x_{\text{Conv}}}_{2}}$) for various values of angular frequency $\omega$.
    No degradation of the solution of FDS was observed in the tested range of angular frequencies.
  }
  \label{fig:omega_error}
\end{figure}
First, $\omega$ was varied from $0.5$ to $14.5$ in increments of $0.5$.
$N$ was set to $6400$.
Figure \ref{fig:omega_time} and Figure \ref{fig:omega_error} respectively show the elapsed time (in seconds) and relative error of the solution vector of the proposed FDS for various angular frequencies.
In Figure \ref{fig:omega_time}, the elapsed times for FDS and Conv increase with increasing angular frequency; however, the increase for FDS is slower than that for Conv.
The reason for this increase in computation time is that it takes longer to compute the Hankel function required in the discretized layer potentials at larger angular frequencies.
Figure \ref{fig:omega_error} shows that the accuracy of the solution computed by FDS is maintained even as the angular frequency increases, at least in the low-frequency range.
Figure \ref{fig:omega_time} and Figure \ref{fig:omega_error} confirm the validity and effectiveness of the proposed fast direct solver.

\begin{figure}[tb]
  \centering
  \includegraphics[width=1.0\linewidth]{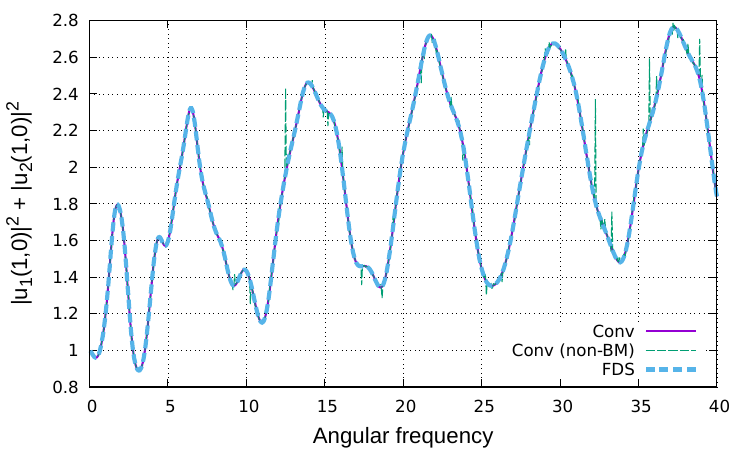}
  \caption{Wave intensity $|u_{1}(x)|^2 + |u_{2}(x)|^2$ at $x = (1, 0)$ for $\omega = 0.05$ to $40$ computed by Conv, Conv (non--BM), and FDS.
    $\omega$ was increased in increments of 0.05.
    Only Conv (non--BM) has spikes, which result from fictitious frequencies.
    FDS and Conv have no spikes (i.e., no fictitious frequencies on real $\omega$).
  }
  \label{fig:wave_intensity}
\end{figure}
Finally, we examine the effect of fictitious eigenvalues on solution accuracy.
Figure \ref{fig:wave_intensity} shows the wave intensity $|u_{1}(x)|^2 + |u_{2}(x)|^2$ at $x = (1, 0)$ computed by Conv, FDS, and a non--Burton--Miller formulated BEM, referred to as Conv (non--BM), respectively.
For Conv (non--BM), the boundary integral equation is expressed as
\begin{align}
  \qty(\mathcal{D}_{ij} + \frac{\delta_{ij}}{2}) u_{j} (x) = u^{I}_{i}(x), \quad x \in \Gamma \label{eq:bie_non_bie}.
\end{align}
In Conv (non--BM), a discretized linear equation of above equation was solved with a non-accelerated LU decomposition.
Figure \ref{fig:wave_intensity} shows that FDS is free from fictitious eigenfrequencies on real $\omega$ (i.e., the wave intensity computed by FDS has no spikes). In contrast, 
Conv (non--BM) has many spikes due to fictitious eigenfrequencies.

\section{Conclusion} \label{sec:conc}
This paper proposed an $O(N)$ fast direct solver for 2D elastic wave scattering problems.
The proposed method is a variant of the method proposed by Martinsson and Rokhlin \cite{martinsson2005fast}.
Our method does not explicitly compute the inverse of the coefficient matrix; therefore, the number of matrix--matrix multiplications required is reduced.
Moreover, the proxy surface method is extended to elastodynamics to obtain shared coefficients for low-rank approximations from discretized integral operators.
Numerical examples show that the proposed method can solve multiple right-hand sides very efficiently, which is a feature of fast direct methods, that it has a complexity of $O(N)$ in the low-frequency range, and that it is capable of very efficient parallel computation.
The proposed method, in combination with the Galerkin discretization of the Barton--Miller--type boundary integral equation, is a versatile fast direct solver for elastic wave scattering problems without fictitious eigenfrequencies in real frequency analysis.

An extension to three dimensions and application to more complicated boundary conditions (e.g., transmission problems \cite{10.1063/1.523808}) will be considered in future research.
To extend the proposed method to three-dimensional problems, we need to change the solver so as to achieve linear time complexity in three dimensions \cite{sushnikova2023fmm}.
In addition, as pointed out in the numerical examples (Section \ref{sec:error}), it is also desirable to improve the accuracy of the multi-level algorithm by modifying the skeleton selection scheme in interpolative decomposition.

\section*{CRediT authorship contribution statement}
{\bf Yasuhiro Matsumoto}:
Conceptualization,
Funding acquisition,
Investigation,
Methodology,
Project administration,
Software,
Resources,
Supervision,
Validation,
Visualization,
Writing -- original draft,
Writing -- review \& editing.
{\bf Taizo Maruyama}:
Funding acquisition,
Methodology,
Software,
Writing -- review \& editing.

\section*{Acknowledgements}
This work was supported by the Japan Society for the Promotion of Science KAKENHI (grant numbers 23K19972 and 24K20783).
This work used the computational resources of the supercomputer TSUBAME4.0 provided by the Institute of Science Tokyo and the supercomputer Camphor3 provided by the Academic Center for Computing and Media Studies of Kyoto University through the Joint Usage/Research Center for Interdisciplinary Large-scale Information Infrastructures and High Performance Computing Infrastructure in Japan (Project IDs: jh240031 and jh240035).

\end{document}